\documentclass[12pt]{amsart}

\def \fg{{\mathfrak g}}

\def \ft{{\mathfrak t}}

\def \u{{\mathfrak u}}

\def \C{{\mathbb C}}

\def \R{{\mathbb R}}

\newtheorem{theorem}{Theorem}[section]
\newtheorem{lemma}[theorem]{Lemma}
\newtheorem{proposition}[theorem]{Proposition}

\numberwithin{equation}{section}

\begin{document}

\baselineskip=16pt

\title[Eigenvalues of sums of pseudo-Hermitian matrices]{Eigenvalues of sums of pseudo-Hermitian matrices}

\author[P. Foth]{Philip Foth}

\address{Department of Mathematics, University of Arizona, Tucson, AZ 85721-0089}

\email{foth@math.arizona.edu}

\subjclass{Primary 15A42, secondary 53D20.}

\keywords{Eigenvalue, pseudo-Hermitian, admissible, convexity.}

\date{May 03, 2008}

\begin{abstract}
 We study analogues of classical inequalities for the eigenvalues of sums of pseudo-Hermitian matrices.
\end{abstract}

\maketitle

\section{Introduction} 

The classical triangle inequality says that for a triangle with side lengths $a$, $b$ and $c$, one has 
$|a-b|\le c\le |a+b|$. If one considers the space $\R^3$ with the Minkowski norm 
$|(x,y,z)|^2=z^2-x^2-y^2$, then in the future timelike cone, defined by $z^2-x^2-y^2 >0 $, $z>0$, the triangle inequality gets reversed, and the sides of a triangle 
${\overrightarrow a}+{\overrightarrow b}={\overrightarrow c}$ satisfy 
$|{\overrightarrow c}| \ge |{\overrightarrow a}|+|{\overrightarrow b}|$. 
This can be interpreted in terms of $2\times 2$ 
traceless pseudo-Hermitian matrices, if one 
puts into correspondence to a vector with coordinates $(x,y,z)$ the matrix 
$$
\left( \begin{array}{cc} z & x+\sqrt{-1}\cdot y \\ -x+\sqrt{-1}\cdot y & -z \end{array} \right) \ .
$$
The eigenvalues of this matrix are $\pm\sqrt{z^2-x^2-y^2}$ and therefore the Minkowski triangle inequality answers 
the following question: given two traceless pseudo-Hermitian matrices with real spectra $(a, -a)$ and $(b, -b)$
and non-negative upper-left entries, what are the possible eigenvalues of their sum? 
Explorations of this and related questions
for Hermitian symmetric matrices (and more generally for triangles in dual vector spaces of compact Lie algebras) 
led to many exciting developments bridging across algebra, Lie theory, representation theory, symplectic geometry, geometric invariant theory, vector bundles, and combinatorics, see for example, \cite{Ful}, \cite{KLM} and references therein. A brief answer to this question can be formulated as follows: given two 
Hermitian symmetric matrices $A$ and $B$,
the set of eigenvalues for their sum $A + B$ necessarily belongs to a convex polytope defined by certain linear inequalities on the sets of eigenvalues of $A$ and $B$. 

In the present paper, we begin answering a similar question in the non-compact setting. Let $G={\rm U}(p,q)$ 
be the pseudounitary Lie group, $\fg$ its Lie algebra, and $\fg^*$ the dual vector space identified with 
the space of pseudo-Hermitian matrices $A$, defined by the condition $A=J_{pq}A^*J_{pq}$, where 
$J_{pq}={\rm diag}(\underbrace{1, ..., 1}_p, \underbrace{-1, ..., -1}_q)$ and $A^*$ is the conjugate transpose.
In general, eigenvalues of pseudo-Hermitian matrices are not necessarily real, unless $A$ is elliptic. 
And moreover, the eigenvalues of the sum of even two elliptic elements can be pretty much arbitrary complex numbers. 
However, if one restricts to the convex cone of {\it admissible} elements \cite{Neeb}, then the question about
possible eigenvalues of the sum becomes more meaningful. In our situation, the convex cone of admissible elements
$\fg^*_{\rm adm}$ will consist of matrices, which are $G$-conjugate to diagonal (and thus real) matrices 
${\rm diag}(\lambda_1,...,\lambda_p, \mu_1, ..., \mu_q)$ such that $\lambda_i > \mu_j$ for all pairs $i,j$. We can 
certainly assume that $\lambda$'s are arranged in the non-increasing order 
$\lambda_1\le\lambda_2\le \cdots \le \lambda_p$ and $\mu$'s are in the non-decreasing order 
$\mu_1\ge \mu_2\ge\cdots\ge\mu_q$ (this is done for convenience), and thus the condition of admissibility becomes rather simple: $\lambda_1 > \mu_1$. 

For two admissible matrices $A, B\in\fg^*_{\rm adm}$ with given spectra, the question of finding
possible eigenvalues of their sum can be formulated in terms of the non-abelian convexity theorem in symplectic
geometry. The coadjoint orbits ${\mathcal O}_A$ and ${\mathcal O}_B$ of $A$ and $B$ carry natural 
invariant symplectic structures and so does their product ${\mathcal O}_A\times {\mathcal O}_B$. 
A generalization due to Weinstein \cite{Wein} of the original Kirwan's theorem to the case of 
non-compact semisimple groups implies that the possible spectrum of $A+B$ forms a convex polyhedral set 
in the positive Weyl chamber $\ft^*_+$ of the dual space to the diagonal torus.   

The primary purpose of this note is to reveal some of the defining conditions on this set, in particulr obtaining 
an analogue of classical Lidskii-Wielandt inequalities \cite{Wi}. Let us formulate our result 
and explain its geometric meaning. For $A, B\in\fg^*_{\rm adm}$ and $C=A+B$, let $\lambda_i(A)$, $\mu_j(A)$,
$\lambda_i(B)$, $\mu_j(B)$, $\lambda_i(C)$, $\mu_j(C)$ be their eigenvalues in the order as above. 
Then for each $m$ integers $1\le i_1< i_2 < \cdots < i_m\le p$ and $\ell$ integers 
$1\le j_1 < j_2 <\cdots < j_\ell \le q$ we have 
$$
\sum_{k=1}^m \lambda_{i_k}(C) \ge \sum_{k=1}^m \lambda_{i_k}(A) + \sum_{k=1}^m \lambda_k(B)
$$
and 
$$
\sum_{k=1}^\ell \mu_{j_k}(C) \le \sum_{k=1}^\ell \mu_{j_k}(A) + \sum_{k=1}^\ell \mu_k(B)\ .
$$
Of course, in addition, we have the trace condition:
$$
\sum_{i=1}^p \lambda_{i}(C)+\sum_{j=1}^q \mu_{j}(C) = 
\sum_{i=1}^p \lambda_{i}(A)+\sum_{j=1}^q \mu_{j}(A) + 
\sum_{i=1}^p \lambda_{i}(B)+\sum_{j=1}^q \mu_{j}(B) \ .
$$
We also state a more general analogue of Thompson-Freede inequalities \cite{TF}. 
Recall from \cite[Theorem VIII.1.19]{Neeb} that the set of possible diagonal entries of an admissible matrix 
$A$ with eigenvalues $({\overrightarrow \lambda}, {\overrightarrow \mu})$ as above, 
form a convex polyhedral set ${\mathcal S}_A$, 
which can be described as the sum $\Pi+{\mathcal C}$ of a polytope $\Pi$ and a cone ${\mathcal C}$. 
The polytope $\Pi$ is the convex hull of 
$S_p\times S_q.({\overrightarrow \lambda}, {\overrightarrow \mu})$ - so its vertices are
obtained by the action of the Weyl group for the maximal compact subgroup (the product of two symmetric groups
in our case). The cone ${\mathcal C}$ is given by the non-compact roots, which in our case means that it is 
the $\R_+$-span of the diagonal differences $a_{ii}-a_{jj}$ for $1\le i\le p$ and $p+1\le j\le n$. 
The above inequalities have then the following geometric interpretation: possible eigenvalues of $A+B$ belong 
to the convex polyhedral region $({\overrightarrow \lambda}(A), {\overrightarrow \mu}(A))+ {\mathcal S}_B$ 
(of course, due to symmetry, we can interchange $A$ and $B$ and get another set of conditions).  

In this note we only deal with analogues of classical eigenvalue inequalities, leaving out natural questions of  relationship with tensor products of representations of $G$ and combinatorics. 

\section{Courant-Fischer theorem for pseudo-Hermitian matrices}

Let $p$ and $q$ be non-negative integers, $n=p+q$, 
and let $G={\rm U}(p,q)$ be the pseudounitary group of $n\times n$ matrices $M$,
satisfying $MJ_{pq}M^*=J_{pq}$, where $J_{pq}$ is the diagonal matrix
$$
J_{pq}=\left( \begin{array}{cc} 1_p & 0 \\ 0 & -1_q \end{array} \right) \ .
$$
Let $\fg=\u(p,q)$ be its Lie algebra of matrices $B$, satisfying $BJ_{pq} + J_{pq}B^*=0$ and let 
$\fg^*$ be its dual vector space, which is identified with the space $\sqrt{-1}\cdot\fg$ of 
pseudo-Hermitian matrices $A$, satisfying $AJ_{pq}=J_{pq}A^*$. In the block form, 
$$
A=\left(
\begin{array}{cc}
H_p & B \\ 
{} & {} \\
-{\bar B^T} & H_q
\end{array}
 \right)\ \ ,
$$
where $H_p$ and $H_q$ are $p\times p$ and $q\times q$ Hermitian symmetric matrices 
respectively and $B$ is a complex $p\times q$ matrix. Let $\fg^*_{\rm adm}$ denote a 
convex component of the open cone of admissible elements, in the terminology of \cite{Neeb}. In general, an element $A\in\fg^*$ is said to be admissible if the co-adjoint orbit ${\mathcal O}_A$ is closed and its convex hull contains no lines. In the pseudounitary case, this translates to the requirement that the coadjoint orbit of $A$ contains a diagonal matrix 
$\Lambda={\rm diag}(\lambda_p, ..., \lambda_1, \mu_1, .., \mu_q)$, where 
$\lambda_p\ge\lambda_{p-1}\ge\cdots\ge\lambda_1$, $\mu_1\ge\mu_2\ge\cdots\ge\mu_q$, and either $\lambda_1 > \mu_1$, or 
$\mu_q > \lambda_p$. There are two open cone components, and without loss of generality we choose $\fg^*_{\rm adm}$ to
be the component in which $\lambda_1 > \mu_1$. 

Let us consider the complex vector space $\C^n$ with the pseudo-Hermitian pairing of signature $(p,q)$:
$$
\langle{\bf z}, {\bf w}\rangle = \sum_{i=1}^p z_i{\bar w}_i - \sum_{j=p+1}^n z_j{\bar w}_j \ .
$$
If we introduce the notation 
$$
{\bf x}^\dagger = (J_{pq}{\bar{\bf x}})^T \ , 
$$
then we can rewrite the above pairing in terms of the usual product:
$$
\langle{\bf z}, {\bf w}\rangle = {\bf w}^\dagger\cdot{\bf z} \ .
$$
Let us also denote by $\C^n_+$ the open cone of positive vectors, satisfying $\langle {\bf z}, {\bf z}\rangle >0$, and
similarly by $\C^n_-$ the cone of negative vectors. Our condition that $A$ is admissible is equivalent to saying that it has real eigenvalues, and the $p$ eigenvalues corresponding to the eigenvectors in $\C^n_+$ are larger than the $q$ eigenvalues corresponding to the eigenvectors in $\C^n_-$. 

Now we shall examine an appropriate analogue of the Rayleigh-Ritz ratio, defined as
$$
{\mathcal R}_A({\bf x})=\frac{{\bf x}^\dagger A{\bf x}}{{\bf x}^\dagger{\bf x}} \ .
$$

\begin{lemma}Let $A\in\fg^*_{\rm adm}$ have the eigenvalues 
\begin{equation}
\lambda_p\ge\lambda_{p-1}\ge\cdots\ge\lambda_1 > \mu_1 \ge \mu_2\ge\cdots\ge \mu_q\ .
\label{eq:eigens}
\end{equation}
Then one has 
$$
\lambda_1=\min_{{\bf x}\in\C^n_+}{\mathcal R}_A({\bf x}) \ \ {\rm and} \ \ 
\mu_1=\max_{{\bf x}\in\C^n_-}{\mathcal R}_A({\bf x}) \ .
$$ 
\label{lem:l2144}
\end{lemma}

\noindent{\it Proof.} Let $U\in G$ be such a matrix that $A=U\Lambda U^{-1}$, where $\Lambda$, as before, is the 
diagonal matrix $\Lambda={\rm diag}(\lambda_p, ..., \lambda_1, \mu_1, .., \mu_q)$. Note that $U^{-1}=U^\dagger$, where
$U^{\dagger}=J_{pq}U^*J_{pq}$. Since $U\in G$, the group of linear transformations of $\C^n$, preserving the pairing $\langle{\bf z},{\bf w}\rangle$, its action on $\C^n$ preserves $\C^n_+$ and $\C^n_-$. For ${\bf x}\in\C^n_+$, denote 
${\bf y}=U^\dagger{\bf x}$, ${\bf y}\in\C^n_+$. Since ${\bf x}^\dagger{\bf x}={\bf y}^\dagger{\bf y} >0$ and 
$(U^\dagger{\bf x})^\dagger = {\bf x}^\dagger U$, we have
$$
{\mathcal R}_A({\bf x})=\frac{{\bf x}^\dagger A{\bf x}}{{\bf x}^\dagger{\bf x}}=\frac{{\bf y}^\dagger\Lambda{\bf y}}{{\bf y}^\dagger{\bf y}}.
$$
Then we need to show that 
$$
{\bf y}^\dagger\Lambda{\bf y}\ge \lambda_1{\bf y}^\dagger{\bf y},
$$
which trivially follows from (\ref{eq:eigens}).

The second statement for $\mu_1$ follows from the statement for $\lambda_1$, by changing $A$ to $-A$. {\bf Q.E.D.} \medskip

Next, let ${\bf v}_1$, .., ${\bf v}_p$, ${\bf w}_1$, ..., ${\bf w}_q$ be a basis of eigenvectors of $A$ in $\C^n$, 
corresponidng to the eigenvalues $\lambda_1$, ..., $\lambda_p$, $\mu_1$, ..., $\mu_q$ respectively and 
orthonormal with respect to $\langle\cdot,\cdot\rangle$. In particular, we have that $||{\bf v}_i||^2=1$, 
$||{\bf w}_j||^2=-1$ and the pairing of any two different vectors from this basis equals zero. Let also, for convenience, denote $V={\rm Span}\{ {\bf v}_1, ..., {\bf v}_p\}$ and  $W={\rm Span}\{ {\bf w}_1, ..., {\bf w}_q\}$.
Note that for $$
{\bf x}=\alpha_1{\bf v}_1 +\cdots + \alpha_p {\bf v}_p 
+\beta_1{\bf w}_1 + \cdots \beta_q{\bf w}_q
$$
the quotient ${\mathcal R}_A({\bf x})$ can be written as
$$
{\mathcal R}_A({\bf x}) = \frac{{\bf x}^\dagger A{\bf x}}{{\bf x}^\dagger{\bf x}}
= \frac{\sum_{i=1}^p |\alpha_i|^2\lambda_i - \sum_{j=1}^q|\beta_j|^2\mu_j}
{\sum_{i=1}^p |\alpha_i|^2 - \sum_{j=1}^q|\beta_j|^2} \ .
$$

From the previous Lemma and the fact that 
$\langle\cdot,\cdot\rangle$ restricts to a positive definite Hermitian pairing on the subspace $V$, which is 
orthogonal to $W$ with respect to $\langle\cdot,\cdot\rangle$, we deduce:

\begin{lemma} 
\begin{equation}
\lambda_k=\min_{{\bf x}\in\C^n_+, \ {\bf x}\perp{\bf v}_1, ..., {\bf v}_{k-1}}{\mathcal R}_A({\bf x}) \ \ \
{\rm and} \ \ \
\lambda_k=\max_{{\bf x}\in V\setminus \{ 0\},\ {\bf x}\perp{\bf v}_{k+1}, ..., {\bf v}_{p}}{\mathcal R}_A({\bf x})\ .
\label{eq:poq}
\end{equation}
\end{lemma}

A similar statement is, of course, valid for $\mu_k$'s: 
$$
\mu_k=\max_{{\bf x}\in\C^n_-, \ {\bf x}\perp{\bf w}_1, ..., {\bf w}_{k-1}}{\mathcal R}_A({\bf x}) \ \ \
{\rm and} \ \ \
\mu_k=\min_{{\bf x}\in W\setminus \{ 0\},\ {\bf x}\perp{\bf w}_{k+1}, ..., {\bf w}_{q}}{\mathcal R}_A({\bf x})\ .
$$

Now we are ready to state and prove a result, similar to the classical Courant-Fischer theorem.

\begin{theorem} Let $A\in\fg^*_{adm}$ be an admissible pseudo-Hermitian matrix with eigenvalues as in (\ref{eq:eigens}).
Let $k$ be an integer, $1\le k\le p$. Then 
\begin{equation}
\lambda_k = \min_{{\bf u}_1, ..., {\bf u}_{n-k}\in\C^n} \ \ \ \ 
\max_{{\bf x}\in\C^n_+, \ {\bf x}\perp{\bf u}_1, ..., {\bf u}_{n-k}} \ \ {\mathcal R}_A({\bf x})
\label{eq:cf1}
\end{equation}
\begin{equation}
\lambda_k = \max_{{\bf u}_1, ..., {\bf u}_{k-1}\in\C^n} \ \ \ \ 
\min_{{\bf x}\in\C^n_+, \ {\bf x}\perp{\bf u}_1, ..., {\bf u}_{k-1}}\ \ {\mathcal R}_A({\bf x})
\label{eq:cf2}
\end{equation}
\label{th:cf}
\end{theorem}

\noindent{\it Proof.} Our line of proof follows the standard argument for the classical Courant-Fischer theorem \cite{HJ}. We will only consider (\ref{eq:cf1}), as the second equality is similar. As in Lemma \ref{lem:l2144}, 
let ${\bf y}=U^\dagger{\bf x}$, where $A=U\Lambda U^{-1}$, and $\Lambda={\rm diag}(\lambda_p, ..., \lambda_1, \mu_1, .., \mu_q)$. Then 
$$
\sup_{{\bf x}\in\C^n_+,\ {\bf x}\perp{\bf u}_1, ..., {\bf u}_{n-k}} \ \ {\mathcal R}_A({\bf x}) \ = \ 
\sup_{{\bf y}\in\C^n_+,\ {\bf y}\perp U^\dagger{\bf u}_1, ..., U^\dagger{\bf u}_{n-k}} \ \ 
{\mathcal R}_{\Lambda}({\bf y}) 
$$
$$
\ge \ \sup_{{\bf y}\in\C^n_+,\ {\bf y}\perp U^\dagger{\bf u}_1, ..., U^\dagger{\bf u}_{n-k},\  
y_{p-k+1}=\cdots=y_{p}=0} \ \ {\mathcal R}_{\Lambda}({\bf y}) \ \ge \ \lambda_k\ .
$$
But (\ref{eq:poq}) shows that the equality holds if we take ${\bf u}_i = {\bf w}_i$ for $1\le i \le q$ and 
${\bf u}_i = {\bf v}_{k-q+i}$ for $q+1\le i \le n-k$. Thus 
$$
\lambda_k = \min_{{\bf u}_1, ..., {\bf u}_{n-k}\in\C^n} \ \ \ \ 
\sup_{{\bf x}\in\C^n_+, \ {\bf x}\perp{\bf u}_1, ..., {\bf u}_{n-k}} \ \ {\mathcal R}_A({\bf x}),
$$
and (\ref{eq:cf2}) is similar. \ \ {\bf Q.E.D.}\medskip

Note that, in general, the ratio ${\mathcal R}_A({\bf x})$ is not bounded from above on $\C^n_+$. Therefore in the
right hand side of the formula (\ref{eq:cf1}), the maximum should be taken over the $(n-k)$-tuples of vectors for which it is actually achieved, and otherwise one might want to use $\sup$ instead of $\max$. 

The above theorem obviously has a natural counterpart, consisting of two series of minimax and maximin identities, for 
$\mu_k$'s. We omit stating and proving those, since it can easily be done if one replaces $A$ by its negative. 

It is also worth noticing that one can rewrite the equality (\ref{eq:cf1}) in the following form:
\begin{equation}
\lambda_k = \ \min_{W_k} \ 
\max_{{\bf x}\in\C^n_+,\ {\bf x}\in W_k} \ {\mathcal R}_A({\bf x}) \ ,
\label{eq:e62}
\end{equation}
where $W_k$ is a subspace of dimension $k$, which in fact can be taken entirely lying in $\C^n_+$ 
(with the exception of the origin, of course). 
 
Next, we state a result similar to one found in \cite{Fan}. 
We will omit the proof since it is a repetition of a standard argument:

\begin{proposition}
For an admissible pseudo-Hermitian matrix $A$ as above, and a positive integer $k\le p$, one has
$$
\lambda_1+\lambda_2+\cdots+\lambda_k = \ \min_{\langle {\bf x}_i, {\bf x}_j\rangle=\delta_{ij}}\ 
\sum_{i=1}^k{\mathcal R}_A({\bf x}_i)\ .
$$
\label{prop:p771}
\end{proposition}

Note that the condition $\langle{\bf x}_i, {\bf x}_j\rangle=\delta_{ij}$ automatically implies that all
of the ${\bf x}_i$'s belong to $\C^n_+$.

As another easy corollary to Theorem \ref{th:cf}, we have the following analogue of classical Weyl inequalities:

\begin{proposition} Let $A, B \in\fg^*_{\rm adm}$ and let $\lambda_i(A)$, $\mu_j(A)$, $\lambda_i(B)$, $\mu_j(B)$, $\lambda_i(A+B)$, $\mu_j(A+B)$ be the eigenvalues of $A$, $B$, and $A+B$ arranged in the order as in 
(\ref{eq:eigens}). Then for each $1\le k\le p$ and $1\le \ell\le q$ we have:
$$
\lambda_k(A+B)\ge \lambda_k(A) + \lambda_1(B) \ \ \ {\rm and} \ \ \ 
\mu_\ell(A+B)\le \mu_\ell(A) + \mu_1(B)\ .
$$
\end{proposition}

\noindent{\it Proof.} We will only prove the first inequality, as the second is similar. 
We know that for each ${\bf x}\in\C^n_+$, one has ${\mathcal R}_B({\bf x})\ge \lambda_1(B)$. 
Hence, using the linearity property of the ratio ${\mathcal R}_{A+B}({\bf x})={\mathcal R}_{A}({\bf x})
+ {\mathcal R}_{B}({\bf x})$,
for $1\le k\le p$ we have 
$$
\lambda_k(A+B) = \ \min_{{\bf u}_1, ..., {\bf u}_{n-k}\in\C^n} \ \ \ \ 
\max_{{\bf x}\in\C^n_+, \ {\bf x}\perp{\bf u}_1, ..., {\bf u}_{n-k}} \ \ {\mathcal R}_{A+B}({\bf x})
$$
$$
\ge\ \min_{{\bf u}_1, ..., {\bf u}_{n-k}\in\C^n} \ \ \ \ 
\max_{{\bf x}\in\C^n_+, \ {\bf x}\perp{\bf u}_1, ..., {\bf u}_{n-k}} \ \ \left( {\mathcal R}_{A}({\bf x})+ \lambda_1(B)\right) \ =\lambda_k(A)+\lambda_1(B)\ .
$$
{\bf Q.E.D.} \medskip

\section{Lidskii-Wieland and Thompson-Freede type inequalities}

In this section we will establish stronger inequalities for the eigenvalues of the sum of two 
admissible pseudo-Hermitian matrices. The first goal of this section is to prove the following
\begin{theorem}
Let $A, B \in\fg^*_{\rm adm}$ and let $\lambda_i(A)$, $\mu_j(A)$, $\lambda_i(B)$, $\mu_j(B)$, $\lambda_i(C)$, $\mu_j(C)$ be the eigenvalues of $A$, $B$, and $C=A+B$ arranged in the order as in 
(\ref{eq:eigens}). Then for each $m$ integers $1\le i_1< i_2 < \cdots < i_m\le p$ and $\ell$ integers 
$1\le j_1 < j_2 <\cdots < j_\ell \le q$ we have 
\begin{equation}
\sum_{k=1}^m \lambda_{i_k}(C) \ge \sum_{k=1}^m \lambda_{i_k}(A) + \sum_{k=1}^m \lambda_k(B)
\label{eq:lw1}
\end{equation}
and 
\begin{equation}
\sum_{k=1}^\ell \mu_{j_k}(C) \le \sum_{k=1}^\ell \mu_{j_k}(A) + \sum_{k=1}^\ell \mu_k(B)\ .
\label{eq:lw2}
\end{equation}
\label{th:lw}
\end{theorem}

In what follows, we will only work on proving (\ref{eq:lw1}), as (\ref{eq:lw1}) is similar. 

For $m\le p$,  us have a fixed $m$-tuple of integers $1\le i_1 < i_2 < \cdots < i_m \le p$. Consider a flag of subspaces $V_{i_1}\subset V_{i_2}\subset \cdots\subset V_{i_m}$, where $V_{i_j}\setminus\{ 0\}\subset \C^n_+$ and
the subscript indicates the dimension of the corresponding subspace. We say that an orthogonal set of vectors
$\{ {\bf x}_{i_1}$, ${\bf x}_{i_2}$, ..., ${\bf x}_{i_m}\}$ is {\it subordinate} to this flag, if 
${\bf x}_{i_j}\in V_{i_j}$ and $\langle {\bf x}_{i_j}, {\bf x}_{i_{k}}\rangle = \delta_{jk}$. 

Denote by $P_m$ the projection operator onto the $Y={\rm Span}\{ {\bf x}_{i_1}, {\bf x}_{i_2}, ..., {\bf x}_{i_m}\}$. 
Here the projection is taken with respect to $\langle\cdot,\cdot\rangle$, and is therefore given by the matrix 
${\bf X}{\bf X}^\dagger$, where the $j$-the column of ${\bf X}$ is ${\bf x}_{i_j}$. For any $A\in\fg^*$, 
the operator $P_mAP_m$ is also pseudo-Hermitian, but its restriction to $Y$ is actually Hermitian, and we let 
$\eta_1\le \eta_2\le\cdots\le \eta_m$ denote the set of its eigenvalues. We have the following analogue of a classical result of Wielandt \cite{Wi}:

\begin{lemma} For $A\in\fg^*_{\rm adm}$ with eigenvalues as in (\ref{eq:eigens}), and $\eta_i$'s as above, we have
$$
\sum_{j=1}^m \lambda_{i_j} =\ \min_{V_{i_1}\subset V_{i_2}\subset \cdots\subset V_{i_m}} \
\max_{{\bf x}_{i_j}\in V_{i_j}} \ \sum_{j=1}^m \eta_j\ .
$$
\label{lemma:wi}
\end{lemma}  

We postpone proving this rather technical lemma till the next secion, and now state an easy corollary:

\begin{proposition} For $A\in\fg^*_{\rm adm}$ with eigenvalues as in (\ref{eq:eigens}) and an $m$-tuple
of integers $1\le i_1 < i_2 < \cdots < i_m \le p$, one has 
\begin{equation}
\sum_{j=1}^m \lambda_{i_j} =\ \min_{V_{i_1}\subset V_{i_2}\subset \cdots\subset V_{i_m}} \
\max_{{\bf x}_{i_j}\in V_{i_j}} \ \sum_{j=1}^m {\mathcal R}_A({\bf x}_{i_j})\ .
\label{eq:e87659}
\end{equation}
\label{prop:p4}
\end{proposition}

\noindent{\it Proof.} One can easily see that the right-hand side of (\ref{eq:e87659}) is exactly the trace of the 
Hermitian operator $P_mAP_m$ acting on the space $Y$, because
$$
\langle P_mAP_m{\bf x}_{i_j}, {\bf x}_{i_k}\rangle = \langle A{\bf x}_{i_j}, {\bf x}_{i_k}\rangle \ ,
$$
and as such, equals $\sum_{j=1}^m \eta_j$. \ {\bf Q.E.D.}\medskip

Now we can establish an analogue of Lidskii-Wieland inequalities.

\noindent{\it Proof of Theorem \ref{th:lw}.} For a given $m$-tuple of integers $1\le i_1< i_2 < \cdots < i_m\le p$, 
let us choose a flag of subspaces $V_{i_1}\subset V_{i_2}\subset \cdots\subset V_{i_m}$ in $\C^n_+$ so that for 
any orthogonal set of vectors
$\{ {\bf x}_{i_1}$, ${\bf x}_{i_2}$, ..., ${\bf x}_{i_m}\}$ subordinate to this flag, one has 
$$
\sum_{j=1}^m \lambda_{i_j}(C) \ge \sum_{j=1}^m {\mathcal R}_C({\bf x}_{i_j}) \ .
$$
As Proposition \ref{prop:p4} shows, this is always possible. Now note that 
$$
\sum_{j=1}^m {\mathcal R}_C({\bf x}_{i_j}) = \sum_{j=1}^m {\mathcal R}_A({\bf x}_{i_j}) + 
\sum_{j=1}^m {\mathcal R}_B({\bf x}_{i_j})\ ,
$$
and use Proposition \ref{prop:p4} once again to choose an orthogonal set of vectors
$\{ {\bf x}_{i_1}$, ${\bf x}_{i_2}$, ..., ${\bf x}_{i_m}\}$ subordinate to the flag
$V_{i_1}\subset V_{i_2}\subset \cdots\subset V_{i_m}$ such that 
$$
\sum_{j=1}^m {\mathcal R}_A({\bf x}_{i_j}) \ge \sum_{j=1}^m \lambda_{i_j}(A).
$$
Next, note that Proposition \ref{prop:p771} implies that 
$$
\sum_{j=1}^m {\mathcal R}_B({\bf x}_{i_j}) \ge \sum_{j=1}^m \lambda_{j}(B) \ ,
$$
and the result follows. {\bf Q.E.D.} \medskip

We now state an analogue of Thompson-Freede inequalities \cite{TF} (without proof). Let us have two $m$-tuples 
of integers $1\le i_1 < i_2 < \cdots < i_m \le p$ and $1\le j_1 < j_2 < \cdots < j_m \le p$ such that 
$i_m+j_m \le m+p$. Then 
$$
\sum_{h=1}^m \lambda_{i_h+j_h-h}(C) \ge 
\sum_{h=1}^m \lambda_{i_h}(A) +
\sum_{h=1}^m \lambda_{j_h}(B)\ .
$$
A similar inequality can be stated for $\mu$'s as well. 

\section{Proof of Lemma \ref{lemma:wi}}

Following the standard path of proving such results as outlined, for example, in the Appendix by B.V. Lidskii to \cite{Gant}, the lemma will follow if we prove the following two statements:\smallskip 

\noindent I. For any flag of subspaces $V_{i_1}\subset V_{i_2}\subset \cdots\subset V_{i_m}$ in $\C^n_+$, there
exist a subordinate set of vectors $\{ {\bf x}_{i_1}$, ${\bf x}_{i_2}$, ..., ${\bf x}_{i_m}\}$, such that
$$
\sum_{j=1}^m \eta_j \ge \sum_{j=1}^m \lambda_{i_j} \ .
$$
\smallskip

\noindent II. There exists a flag $V_{i_1}\subset V_{i_2}\subset \cdots\subset V_{i_m}$ such that for 
any subordinate set of vectors $\{ {\bf x}_{i_1}$, ${\bf x}_{i_2}$, ..., ${\bf x}_{i_m}\}$, one has 
$$
\sum_{j=1}^m \lambda_{i_j} \ge \sum_{j=1}^m \eta_j \ .
$$
\smallskip

We will first prove II. Set
$$
V_{i_j} = {\rm Span}\{ {\bf v}_1, ..., {\bf v}_{i_j}\}\ , 
$$
where ${\bf v}_1$, .., ${\bf v}_p \in\C^n_+$ are eigenvectors of $A$, 
corresponidng to the eigenvalues $\lambda_1$, ..., $\lambda_p$ respectively. 
Note that $V_{i_j}\setminus\{ 0\}\subset\C^n_+$. Let $\{ {\bf x}_{i_1}$, ${\bf x}_{i_2}$, ..., ${\bf x}_{i_m}\}$
be a set of vectors subordinate to the chosen flag, and let $W_\ell$ be an $\ell$-dimensional subspace in their span. 
We know from the classical minimax identities that 
$$
\eta_\ell \le\ \max_{{\bf x}\in W_\ell}\ {\mathcal R}_{P_mAP_m}({\bf x)}\ . 
$$
Note that for ${\bf x}\in W_\ell$, we have 
${\mathcal R}_{P_mAP_m}({\bf x}) = {\mathcal R}_A({\bf x})$. Thus if we let 
$W_\ell = {\rm Span}\{ {\bf x}_{i_1}, {\bf x}_{i_2}, ..., {\bf x}_{i_\ell}\}$, then the fact that 
$W_\ell\subset V_{i_\ell}$ will imply 
$$
\max_{{\bf x}\in W_\ell}\ {\mathcal R}_{A}({\bf x)} \le \max_{{\bf x}\in V_\ell}\ {\mathcal R}_{A}({\bf x)}\ .
$$
But the maximum in the right-hand side is achieved on the eigenvector ${\bf v}_{i_\ell}$ and equals $\lambda_{i_\ell}$. 
(We recall that the operator $A$ is trivially Hermitian on the span of its eigenvectors from $\C^n_+$.) Thus
$$
\eta_\ell \le\ \max_{{\bf x}\in W_\ell}\ {\mathcal R}_{P_mAP_m}({\bf x)} = 
\ \max_{{\bf x}\in W_\ell}\ {\mathcal R}_{A}({\bf x)} \le\ \max_{{\bf x}\in V_\ell}\ {\mathcal R}_{A}({\bf x)} = \lambda_{i_\ell}\ ,
$$
proving II.
\smallskip

Now we turn to proving I, by induction on $p$. Note that for $p=1$, the statement amounts to showing that 
$$
\lambda_1 = \min_{V_1}\ \eta_1\ ,
$$
where $V_1$ is a one-dimensional subspace in $\C^n_+$. This is not hard to establish directly, and in any case, is an easy consequence of \cite[Proposition 4.1]{Fot}. 

Now we can take $m<p$, since in the case when $m=p$, the statement is again a consequence of {\it loc.cit.} We consider two subcases: 

1). When $i_m < p$, there exists a $(p-1)$-dimensional subspace $R_{p-1}$ of $\C^n_+$, containing the
whole flag $V_{i_1}\subset V_{i_2}\subset \cdots\subset V_{i_m}$. Let $P_{p-1}$ be the operator of projection
onto $R_{p-1}$. Consider the pseudo-Hermitian operator $A_{p-1}=P_{p-1}AP_{p-1}$, which is actually Hermitian, being 
restricted to $R_{p-1}$. Clearly for all ${\bf x}\in R_{p-1}$, one has 
${\mathcal R}_{A_{p-1}}({\bf x})={\mathcal R}_A({\bf x})$. If we denote by $\xi_1, .., \xi_{p-1}$ the eigenvalues of 
$A_{p-1}$, in the non-decreasing order, then according to {\it loc.cit.}, one has 
\begin{equation}
\xi_i \ge \lambda_i \ \ \ {\rm for} \ \ 1\le i \le p-1\ .
\label{eq:e54549}
\end{equation}
By the inductive hypothesis, for any flag $V_{i_1}\subset V_{i_2}\subset \cdots\subset V_{i_m}$ in $R_{n-1}$, 
there exists a subordinate system of vectors $\{ {\bf x}_{i_1}$, ${\bf x}_{i_2}$, ..., ${\bf x}_{i_m}\}$ such that 
$$
\sum_{j=1}^m \eta_j \ge \sum_{j=1}^m \xi_{i_j} \ ,
$$
and we are done in this case. 

2). Now consider the case $i_m=p$. Assume $i_m=p$, $i_{m-1}=p-1$, ..., $i_{m-s}=p-s$ and that the number 
$(p-s-1)$ is not a part of the $m$-tuple $1\le i_1 < i_2 < \cdots < i_m \le p$. Let $i_t$ be the largest 
remaining element of this $m$-tuple (the case when there is no such left requires only a minor and 
trivial modification of our discussion). The corresponding flag of subspaces now takes the form
$$
V_{i_1}\subset V_{i_2}\subset\cdots\subset V_{i_t}\subset V_{i_t+1}\subset\cdots\subset V_{p} \ .
$$
Let ${\bf v}_{p-s}$, ${\bf v}_{p-s+1}$, ..., ${\bf v}_p$ be the eigenvectors of $A$ corresponding to the 
$s+1$ largest eigenvalues. Let $R_{n-1}$ be the subspace of $\C^n$ spanned by these vectors and containing 
$V_{i_t}$ and all the ${\bf w}_j$'s. Such a subspace exists since $i_t\le p-s-2$ and thus $s+1+i_t\le p-1$. 

Consider yet another flag of subspaces:
\begin{equation}
V_{i_1}\subset V_{i_2}\subset\cdots\subset V_{i_t}\subset R_{p-s-1}\subset R_{p-s}\subset\cdots\subset
R_{p-1}\ ,
\label{eq:newflag}
\end{equation}
where $R_j=V_{j+1}\cap R_{n-1}$. (In the degenerate case when the dimension of the intersection does not drop by 1, 
we can artificially remove one extra dimension.)

Again, let us introduce the operator $A_{p-1}=P_{p-1}AP_{p-1}$ on the space $R_{p-1}$ as before. Using our inductive assumption, we can find a subordinate system of vectors 
$$
\{ {\bf x}_{i_1}, {\bf x}_{i_2}, ..., {\bf x}_{i_t}, {\bf x}_{p-s-1}, ..., {\bf x}_{p-1}\}
$$
such that 
$$
\sum_{j=1}^m \eta_j \ge \sum_{j=1}^t \xi_{i_j}+\sum_{j=p-s-1}^{p-1}\xi_j\ , 
$$
where $\xi$'s are the eigenvalues of $A_{p-1}$ arranged in the non-decreasing order. 
According to (\ref{eq:e54549}), we have 
$$
\xi_{i_1}\ge \lambda_{i_1},\ \xi_{i_2}\ge \lambda_{i_2},\ ...,\ \xi_{i_t}\ge \lambda_{i_t}\ .
$$
The vectors ${\bf v}_{p-s}$, ${\bf v}_{p-s+1}$, ..., ${\bf v}_p$ belong to the subspace $R_{p-1}$ and
are eigenvectors for $A_{p-1}$. Thus the corresponding eigenvalues $\lambda_{p-s}, ..., \lambda_{p}$ 
are dominated by $\xi_{p-s-1}$, ..., $\xi_{p-1}$, which are the largest $(s+1)$ eigenvalues of $A_{p-1}$. 
Thus we conclude that 
$$
\xi_{i_1}+\xi_{i_2}+\cdots + \xi_{i_t}+\xi_{p-s}+ \cdots + \xi_{p} \ge
\lambda_{i_1}+\lambda_{i_2}+\cdots + \lambda_{i_t}+\lambda_{p-s}+\cdots + \lambda_p
$$
Since the system $\{ {\bf x}_{i_1}, {\bf x}_{i_2}, ..., {\bf x}_{i_t}, {\bf x}_{p-s-1}, ..., {\bf x}_{p-1}\}$
is subordinate not only to the orginal flag, but also to (\ref{eq:newflag}), we have completed the proof.

%\section*{Acknowledgements}
%
%I am grateful to Patrick Sidney for useful discussions.   

%%%%%%%%%%%%%%%%%%%%%%%%%%%%%%%%%%%%%%%%%%%%%%%%%%%%%%%%%%%%%%%%

\end{document}